\documentclass{ws-procs9x6}
\usepackage{pslatex}
\usepackage[ansinew]{inputenc}

\newcommand{\dist}{\operatorname{dist}}

\def \NN {{\mathbb N}}

 \def \be {\beta } 
 \def \vare {\varepsilon }  \def \th {\theta }

  \def \cs {\mathcal{S}}

\newcommand{\dem}{\begin{proof}}
\newcommand{\cqd}{\end{proof}}

\newcommand{\qand}{\quad\text{and}\quad}

\newtheorem{maintheorem}{Theorem}

\newcommand{\cmt}{\begin{maintheorem}}
\newcommand{\fmt}{\end{maintheorem}}

\newtheorem{maincorollary}[maintheorem]{Corollary}
\newcommand{\cmc}{\begin{maincorollary}}
\newcommand{\fmc}{\end{maincorollary}}

\newtheorem{T}{Theorem}[section]
\newcommand{\ct}{\begin{T}}
\newcommand{\ft}{\end{T}}

\newtheorem{Corollary}[T]{Corollary}
\newcommand{\cco}{\begin{Corollary}}
\newcommand{\fco}{\end{Corollary}}

\newtheorem{Proposition}[T]{Proposition}
\newcommand{\cpr}{\begin{Proposition}}
\newcommand{\fpr}{\end{Proposition}}

\newtheorem{Lemma}[T]{Lemma}
\newcommand{\cle}{\begin{Lemma}}
\newcommand{\fle}{\end{Lemma}}

\newtheorem{Remark}[T]{Remark}
\newcommand{\cre}{\begin{Remark}}
\newcommand{\fre}{\end{Remark}}

\newtheorem{Definition}{Definition}
\newcommand{\cde}{\begin{Definition}}
\newcommand{\fde}{\end{Definition}}

\begin{document}

\title{Integrability versus frequency of
hyperbolic times and the existence of a.c.i.m.}



\author{V{\'\i}tor Ara\'ujo}
\address{Centro de Matem\'atica da Universidade do Porto\\
Rua do Campo Alegre 687, 4169-007 Porto, Portugal\\
Email:vdaraujo@fc.up.pt\\
URL:www.fc.up.pt/cmup/home/vdaraujo}

\maketitle

\abstracts{
We consider  dynamical systems on compact manifolds, which are
local diffeomorphisms outside an exceptional set (a compact
submanifold). We are interested in analyzing the relation between
the integrability (with respect to Lebesgue measure) of the first
hyperbolic time map  and the existence of  positive frequency of
hyperbolic times. We show that some (strong) integrability of the
first hyperbolic time map implies   positive frequency of
hyperbolic times. We also present an example of a map with
positive frequency of hyperbolic times at Lebesgue almost every
point but whose first hyperbolic time map is not integrable with
respect to the Lebesgue measure.
}

\section{Introduction and statement of results}
\label{sec:introduction}
This note presents results from a joint work
\cite{alves-araujo2003} with J. F. Alves.

The study of non-uniformly hyperbolic behavior has recently
benefitted from a new tool introduced by
Alves\cite{alves2000} through the notion of {\em hyperbolic
times}.
These have played an
important role in the proof of several results about the
existence of absolutely continuous invariant measures and
their statistical properties; see  Alves\cite{alves2000},
Alves-Araújo\cite{alves-araujo},
Alves-Bonatti-Viana\cite{alves-bonatti-viana2000},
Alves-Luzzatto-Pinheiro\cite{alves-luzzatto-pinheiro}
and Alves-Viana\cite{alves-viana2002}.

Let $f:M\to M$ be a continuous map on
a compact Riemannian manifold $M$ with the induced distance
$\dist$, and fix a normalized Riemannian
volume form $m$ on $M$ that we call {\em Lebesgue measure}.
We assume that \( f \) is a $C^2$ local
diffeomorphism in all of $M$ but an exceptional set
 \( \cs \subset M\), where $\cs$ is a compact submanifold of
 $M$ with $\dim(\cs)<\dim(M)$ (thus $m(\cs)=0$) such that 
$f$  {\em behaves like a power of the distance}
 to \( \cs \): for some $\beta>0$
 \begin{equation}
   \label{eq:nondegenerate}
   \| Df(x) \| \approx \dist(x,\cs)^{-\be}
 \end{equation}
(see Alves-Araújo\cite{alves-araujo2003} for a precise statement).
The set $\cs$ may be the set of critical points of $f$ or a
set where $f$ fails to be differentiable ($\cs=\emptyset$
may also be considered).
Let \(
\dist_{\delta}(x,\cs) =  \dist(x,\cs) \) if \( \dist(x,\cs) \leq
\delta\), and \( \dist_{\delta}(x,\cs) =1 \) otherwise.
\cde\label{def.ht}
 Let  $\beta>0$ be as in~(\ref{eq:nondegenerate}), and
choose $0< b < \min\{1/2,1/(4\beta)\}$.
Given $0<\sigma<1$ and $\delta>0$, we say that $n$ is a {\em
$(\sigma,\delta)$-hyperbolic time\footnote{In the case
$\cs=\emptyset$ the definition of $(\sigma,\delta)$-hyperbolic
time reduces to the first condition in \eqref{d.ht}}}
for $x\in M$ if for every $1\le k \le n$,
 \begin{equation}\label{d.ht}
\prod_{j=n-k}^{n-1}\|Df(f^j(x))^{-1}\| \le \sigma^k \qand
\dist_\delta(f^{n-k}(x), \cs)\ge \sigma^{b k}.
 \end{equation}
We say that the {\em frequency of
$(\sigma,\delta)$-hyperbolic times} for $x\in M$ is greater than
$\theta>0$ if the hyperbolic times 
$n_1<n_2<n_3\dots$ satisfy $\liminf_N N^{-1}\#\{k:n_k < N\} \ge \theta$.
\fde
The following result contained in
Alves-Bonatti-Viana\cite{alves-bonatti-viana2000} is a
strong motivation for the study of
 hyperbolic times.

\cmt\label{t.abv} If there are $0<\sigma<1$,
   $\delta>0$,  and $H\subset M$ with $m(H)>0$ such that
the frequency of $(\sigma,\delta)$-hyperbolic  times is bigger
than $\theta>0$ for every $x\in H$, then $f$ has some absolutely
continuous invariant probability measure.
 \fmt

Assuming $(\sigma,\delta)$-hyperbolic times exist Lebesgue
almost everywhere, let $h:M\to\NN$ be such that $h(x)$ is
the first $(\sigma,\delta)$-hyperbolic time. The
integrability properties of $h$ are important in the study
of stochastic stability and decay of correlations for
several classes of dynamical systems, see
Alves-Araújo\cite{alves-araujo},
Alves-Viana\cite{alves-viana2002} and 
Alves-Luzzatto-Pinheiro\cite{alves-luzzatto-pinheiro}.

 \cmt\label{t.integrable-h-1}
If, for some $0<\sigma<1$ and $\delta>0$, $h$ is
Lebesgue integrable, then $f$ has an  absolutely continuous invariant
probability measure $\mu$.
 \fmt
From these two results it is natural to investigate the
relation between the integrability of $h$ and positive
frequency of hyperbolic times.

\cmt\label{t.integrable-h-2} If for $0<\sigma<1$ and
$\delta>0$ the first $(\sigma,\delta)$-hyperbolic time map
$h$ belongs to $L^p(m)$ for some $p>4$, then there
are $\hat\sigma>0$ and $\theta>0$ such that Lebesgue almost
every $x\in M$ has frequency of
$(\hat\sigma,\delta)$-hyperbolic times bigger than $\theta$.
\fmt

The smallest value of $p\ge 1$ for which the first condition
in Theorem~\ref{t.integrable-h-2} still implies the desired
conclusion is not known (to the best of our knowledge).
However, as a consequence of our proofs we obtain
an optimal bound when $\cs=\emptyset$.

\cmc\label{co:t-integrable-h} Let $f\colon M\to M$ be a $C^2$
local diffeomorphism. If for some $\sigma\in(0,1)$ the first
$\sigma$-hyperbolic time map is Lebesgue integrable, then there
are $\hat\sigma>0$ and $\th>0$   such that Lebesgue almost every
$x\in M$ has frequency of $\hat\sigma$-hyperbolic times bigger
than $\theta$.  \fmc

\begin{figure}[ht]
\centerline{\epsfxsize=4cm\epsfbox{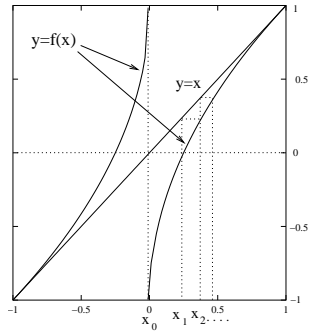}}   
  \caption{A map with non-integrable first hyperbolic
  time map.}
  \label{fig1}
\end{figure}

There is no hope of a result in the oposite direction, at
least when $\cs\neq\emptyset$.
We consider the map $\hat f:I\to I$
(see figure~\ref{fig1}) with $I=[-1,1]$ given by
\[
f(x)= 2\sqrt{x}-1 \quad\text{if $x\ge0$}
\quad \text{and} \quad
f(x)= 1- 2\sqrt{|x|} \quad\text{otherwise}.
\]
This induces a local homeomorphism $f:S^1\to S^1$
not differentiable at the point $0$.

It is not difficult to show that Lebesgue measure is
$f$-invariant and ergodic.  Using this we show that
Lebesgue almost every point has positive frequency of
hyperbolic times. Moreover, it is easy to see that for every
$x$ the first hyperbolic time must be at least the first
return time to a neighborhood of $0$.  Hence $\int h\,dm \ge
\sum_{n\ge1} n(x_{n+1}-x_n)$ where $x_n=g^n(0), n\ge1$ and
\( g:(-1,1)\to(0,1), x\mapsto (1+x)^2/4 \) is an inverse
branch of $f$. It is an easy exercise to check that the
series diverges. Thus $h$ is not Lebesgue integrable (see
Alves-Araujo\cite{alves-araujo2003} for details).

We do not know whether it is possible to construct an
example  with the same properties but in the local
diffeomorphism setting (i.e., with $\cs=\emptyset$).

\section{Sketch of some proofs}
\label{sec:sketch}

We present a sketch of the proofs of Theorems~\ref{t.integrable-h-1}
and~\ref{t.integrable-h-2}.
The strategy is to consider the
weak$^*$ accumulation points of
$\mu_n=n^{-1}\sum_{j=0}^{n-1} f_*^j m$ when
$n\to\infty$,
which are $f$-invariant, and to show that these accumulation
points are absolutely continuous. We use the main feature of
hyperbolic times: their uniformly expanding behavior.

 \cle\label{p.contr}
 There are $\delta_1>0$ and $C_1>0$
 such that if $n$ is a
$(\sigma,\delta)$-hyperbolic time for $x$, then there is a
neighborhood $V_x$ of $x$ in $M$ for which
 \begin{enumerate}
 \item $f^n$ maps $V_x$ diffeomorphically onto
 $B(f^n(x),\delta_1)$;
 \item for $k=1,\dots,n-1$ and $y,z\in V_x$, we have 
\newline
\centerline{$\dist(f^{n-k}(y),f^{n-k}(z)) \le
\sigma^{k/2}\dist(f^{n}(y),f^{n}(z))$;}
 \item $f^n\vert V_x$ has distortion bounded by $C_1$: 
$$\text{if  $y, z\in V_x$, then }
 \frac{1}{C_1} \le \frac{|\det Df^n (y)|}{|\det
Df^n (z)|}\le C_1 \,. $$
  \end{enumerate}
  \fle 
\dem See Lemmas 5.2 and 5.3 of
  Alves-Bonatti-Viana\cite{alves-bonatti-viana2000}. \cqd

From this we obtain the following result in
Alves-Araújo\cite{alves-araujo2003}, Proposition 3.1,
where we write $H_n=\{x\in M: $n$\mbox{ is a
  $(\sigma,\delta)$-hyperbolic time for } x \}$.

\cpr\label{p.density}
 There is a constant $C_2>0$ 
(only depending on $\delta_1$ and $C_1$)
  such that 
 $
 \frac{d}{dm}f^n_*\big(m\mid H_n\big)\leq C_2
 $
for every $n\geq 0$
 \fpr

Defining for each $n$ the sets $H_n^*=h^{-1}(n)$ we have
$m(h)=\sum_{n\ge1} n\cdot m(H_n^*)$. We need to denote
$R_{n,k}=\{x\in H_n: h(f^n(x))=k\}$ and define the measures
$\nu_n=f^n_*(m\mid H_n)$ and
$\eta_n=\sum_{k=2}^\infty\sum_{j=1}^{k-1} f^{n+j}_*(m\mid
R_{n,k})$. Then $\mu_n\le n^{-1}\sum_{j=0}^{n-1}
(\nu_j+\eta_j)$.

From Proposition~\ref{p.density} we have
$\frac{d\nu_n}{dm}\leq C_2$ where $C_2$ does not depend on
$n$.

 \cpr\label{p.dens}
 Given $\vare>0$, there is
$C_3(\vare)>0$ such that for  every $n\geq 1$ there are 
 non-negative measures such that
 $\omega+\rho\ge\eta_n$, with
 $
 \frac{d\omega}{dm}\leq C_3(\vare)$
and
$
 \rho(M)<\vare.
$
 \fpr

\dem We let
$
 \omega=C_2\sum_{k=2}^{\ell-1}\sum_{j=1}^{k-1}
 f^j_*(m\mid H^*_k)
$
and
$
 \rho=C_2\sum_{k=\ell}^\infty\sum_{j=1}^{k-1}
 f^j_*(m\mid H^*_k),
$
where $\ell\in\NN$ is big enough so that
$
 \sum_{j=\ell}^{\infty}k\,
 m\big(H^*_k\big)<C_2^{-1}\vare
$
(by the integrability of $h$).
See Alves-Araújo\cite{alves-araujo2003} for more details.
\cqd

This ensures that every weak$^*$ accumulation point $\mu$ of
$(\mu_n)_n$ must be absolutely continuous, proving
Theorem~\ref{t.integrable-h-1}. 

To prove Theorem~\ref{t.integrable-h-2} we start by a simple
geometric result (Proposition 4.1 in
Alves-Araújo\cite{alves-araujo2003}).

\cpr\label{le:logdist-int} If $\cs$ is a compact submanifold of
$M$ with $\dim(\cs)<\dim(M)$, then the function $\log\dist(x,\cs)$
belongs to $L^p(m)$ for every $1\le p<\infty$. \fpr

The standard H\"older inequality easily implies the following
\cco\label{c.integra} If the density $d\mu/dm$ belongs to $L^q(m)$
for some $q>1$, then $\log\dist(x,\cs)$ is $\mu$-integrable. \fco

\cle\label{le:dist-int} If the first $(\sigma,\delta)$-hyperbolic
time map $h:M\to\NN$ belongs to $L^p(m)$ for some $p>4$, then
$\log\dist(x,\cs)$ is $\mu$-integrable. \fle

\dem
We take $\vare>0$ and use the construction of the
a.c.i.m. $\mu$ and the definitions of $\nu_n, \eta_n$ and
$\rho$.

Let $d_i=\sigma^{bi}$ for $i\ge1$ and $b$ from
Definition~\ref{def.ht}, and define
$ B_i=\{x\in M\colon \dist(x,\cs)\le d_i\}$. Then for $n$ a
$(\sigma,\delta)$-hyperbolic time for $x$ we have
$f^j(x)\in
M\setminus B_i$ for all $j\in\{n-i,\dots,n-1\}$. This
implies for all $i\ge1$
\begin{align*}
\rho(B_i) & =  C_2 \sum_{k=\ell}^\infty\sum_{j=1}^{k-1} m(
H_k^*\cap f^{-j} (B_i ) )
= C_2
\sum_{k=\ell}^\infty\sum_{j=1}^{k-i} m( H_k^*\cap f^{-j} (B_i ))
\\
& \le C_2 \sum_{k=\max\{\ell,i\}}^\infty\sum_{j=1}^{k-i} m(
H_k^*\cap f^{-j} (B_i) ) 
\le C_2  \sum_{k=i}^\infty k \, m(
H^*_k).
\end{align*}
\noindent
Now by Propositions~\ref{p.density} and~\ref{p.dens} we
know that
$
\mu_n\le n^{-1}\sum_{j=0}^{n-1}\nu_j+\omega+\rho\le\nu+\rho
$
where $\nu$ is a measure with uniformly bounded density. Hence any
weak$^*$ accumulation point $\mu$ of  $(\mu_n)_n$ is
bounded by $\nu+\rho$. Since  $\cs$ is a
submanifold of $M$, then $\log\dist(x,\cs)\in L^1(\nu)$
by Proposition~\ref{le:logdist-int}. On the other
hand,
 \begin{align*}
 \int_M -\log\dist_\delta(x,\cs)\,d\rho &\le
\sum_{i=1}^\infty -\rho(B_i)\log d_{i+1}  \le -b \log\sigma
\sum_{i=1}^\infty (i+1) \sum_{k=i}^\infty k \, m( H^*_k).
\end{align*}
We have $h\in L^p(m)$, which means
$\sum_{k\ge1} k^p m(H_k^*)<\infty$. Hence there is
a constant $C>0$ such that $m(H_k^*)\le C k^{-p}$ for all
$k\ge1$. Thus for $i\ge 2$
\[
\sum_{k=i}^\infty k \, m( H^*_k) \le \sum_{k=i}^\infty
\frac{C}{k^{p-1}}\le \int_{i-1}^\infty \frac{C}{x^{p-1}} \, dx =
\frac{C}{(p-2)(i-1)^{p-2}},
\]
and so
\[
\sum_{i=2}^\infty (i+1) \sum_{k=i}^\infty k \, m( H^*_k) \le
\frac{C}{p-2}\sum_{i=2}^\infty\frac{i+1}{(i-1)^{p-2}}.
\]
This last series is summable for $p>4$. Hence
$\log\dist(x,\cs)\in L^1(\mu)$ for $p>4$.
 \cqd

 Since $f$
 behaves like a power of the distance near 
 $\cs$, the following is a corollary of the previous results
 (see Alves-Araújo\cite{alves-araujo2003}, Corollary 4.4).

\cco\label{co:dist-int} If the first $(\sigma,\delta)$-hyperbolic
time map $h$ belongs to $L^p(m)$ for
some $p>4$, 
then $ \log \| Df(x)^{-1} \|$ is $\mu$-integrable. \fco

By definition of hyperbolic time, if $x\in H_n$ and
$f^n(x)\in H_k$, then $x\in H_{n+k}$. Since $f$ preserves
sets of zero $m$ measure and $m$-a.e.
point has a hyperbolic time, then $m$ almost all points must
have infinitely many hyperbolic times. Thus we have
$
\liminf_{n} n^{-1} \sum_{j=0}^{n-1} \log \| Df(f^j(
x))^{-1} \| \le \log\sigma <0
$
for $m$-a.e. $x\in M$, and hence  $\mu$ almost
everywhere.  The $\mu$-integrability of $\log \| Df(x)^{-1}
\|$ and the Ergodic Theorem  ensure that the
limit of the above time average exists  $\mu$-a.e.

If $\cs=\emptyset$, then
$\log \| Df(x)^{-1} \|\in L^1(\mu)$ 
(it is bounded). Hence the negative bound
 on the above time averages
is enough for obtaining Corollary~\ref{co:t-integrable-h} by
applying the following result from
Alves-Bonatti-Viana\cite{alves-bonatti-viana2000}.

\cmt\label{t.abv2} Let $f\colon M\to M$ be a $C^2$ local
diffeomorphism outside an exceptional set
$\cs\subset M$. If  for some set $H\subset M$ with $m(H)>0$
\begin{enumerate}
\item $\exists c>0 \, \forall x\in H \, , \limsup_{n} n^{-1}
\sum_{j=0}^{n-1}\log\|Df(f^j(x))^{-1}\|<-c$;
\item $\forall\vare>0\, \exists \delta>0 \,\forall x\in H,
 \limsup_n
n^{-1} \sum_{j=0}^{n-1}-\log \dist_\delta(f^j(x),\cs)
\le\vare,$
\end{enumerate}
then there are $0<\sigma<1$, $\delta>0$ and $\theta>0$
such that the frequency of $(\sigma,\delta)$-hyperbolic times for
the points in $H$ is bigger than $\theta$. 
 \fmt

An a.c.i.m. $\mu$ is sure to exist from $h\in L^p(m), p>1$,
according to Theorem~\ref{t.integrable-h-1}, and
$\log\dist(x,\cs)\in L^1(\mu)$  after
Corollary~\ref{co:dist-int}.

\cle\label{le:both-averages} 
If $\mu$ is an $f$-invariant probability measure and
$\log\dist(x,\cs)$ is $\mu$-integrable, then 
$\mu$-almost every $x\in M$ satisfies the second condition
of Theorem~\ref{t.abv2}.
\fle 

\dem
Easy when $\mu$ is ergodic. For details see
Alves-Araújo\cite{alves-araujo2003}, Lemma 4.6.
\cqd

Now we just apply Theorem~\ref{t.abv2} to obtain
Theorem~\ref{t.integrable-h-2}. 

\subsection*{Acknowledgments}
 Work partially supported by ESF-PRODYN
 and FCT-CMUP(Portugal).




\begin{thebibliography}{1}

\bibitem{alves2000}
J.~F. Alves.
\newblock {SRB} measures for non-hyperbolic systems with multidimensional
  expansion.
\newblock {\em Ann. Scient. Éc. Norm. Sup.}, 33:1--32, 2000.
\newblock $4^e$ serie.

\bibitem{alves-araujo2003}
J.~F. Alves and V.~Araújo.
\newblock Hyperbolic times: frequency versus integrability.
\newblock {\em Ergodic Theory and Dynamical Systems}, 24:1--18, 2003.

\bibitem{alves-araujo}
J.~F. Alves and V.~Araújo.
\newblock Random perturbations of non-uniformly expanding maps.
\newblock {\em Astérisque}, 286:25--62, 2003.

\bibitem{alves-bonatti-viana2000}
J.~F. Alves, C.~Bonatti, and M.~Viana.
\newblock Srb measures for partially hyperbolic systems whose central direction
  is mostly expanding.
\newblock {\em Invent. Math.}, 140:351--398, 2000.

\bibitem{alves-luzzatto-pinheiro}
J.~F. Alves, S.~Luzzatto, and V.~Pinheiro.
\newblock Markov structures and decay of correlations for non-uniformly
  expanding dynamical systems.
\newblock {\em Preprint CMUP}, 2002.

\bibitem{alves-viana2002}
J.~F. Alves and M.~Viana.
\newblock Statistical stability for robust classes of maps with non-uniform
  expansion.
\newblock {\em Ergodic Theory and Dynamical Systems}, 22:1--32, 2002.

\end{thebibliography}
\end{document}